\documentclass{amsart}
\usepackage{amssymb,latexsym,amsmath,amsfonts,hhline}
\usepackage{pdfsync}
\usepackage{xcolor} 
\usepackage{comment,colortbl}

\newcommand{\cal}{\mathcal}

\newtheorem{formula}{}[section]
\newtheorem{definition}[formula]{Definition}
\newtheorem{corollary}[formula]{Corollary}
\newtheorem{remark}[formula]{Remark}
\newtheorem{lemma}[formula]{Lemma}
\newtheorem{theorem}[formula]{Theorem}
\newtheorem{proposition}[formula]{Proposition}

\def\thrm{\begin{theorem}}
\def\thrml#1{\begin{theorem}\label{#1}}
\def\ethrm{\end{theorem}}
\def\rmrk{\begin{remark}}
\def\rmrkl#1{\begin{remark}\label{#1}}
\def\ermrk{\end{remark}}
\def\dfntn{\begin{definition}}
\def\dfntnl#1{\begin{definition}\label{#1}}
\def\edfntn{\end{definition}}
\def\nmrt{\begin{enumerate}}
\def\enmrt{\end{enumerate}}
\def\tm#1{\item[{\rm (#1)}]}
\def\qtn{\begin{equation}}
\def\qtnl#1{\begin{equation}\label{#1}}
\def\eqtn{\end{equation}}
\def\lmm{\begin{lemma}}
\def\lmml#1{\begin{lemma}\label{#1}}
\def\elmm{\end{lemma}}
\def\crllr{\begin{corollary}}
\def\crllrl#1{\begin{corollary}\label{#1}}
\def\ecrllr{\end{corollary}}
\def\css{\begin{cases}}
\def\ecss{\end{cases}}
\def\prpstn{\begin{proposition}}
\def\prpstnl#1{\begin{proposition}\label{#1}}
\def\eprpstn{\end{proposition}}

\def\proof{\noindent{\bf Proof}.\ }
					
\def\cA{{\cal A}}

\def\cC{{\cal C}}

\def\cS{{\cal S}}

\def\cX{{\cal X}}

\def\mF{{\mathbb F}}

\DeclareMathOperator{\aut}{Aut}
\DeclareMathOperator{\alt}{Alt}
\DeclareMathOperator{\AGL}{AGL}

\DeclareMathOperator{\Aaut}{Aut_{alg}}

\DeclareMathOperator{\GL}{GL}

\DeclareMathOperator{\inv}{Inv}

\DeclareMathOperator{\PGL}{PGL}

\DeclareMathOperator{\sym}{Sym}

\def\eprf{\hfill$\square$}

\def\mmod#1#2#3{#1=#2\ (\text{\rm mod}\hspace{2pt}#3)}
\def\qaq{\quad\text{and}\quad}
\def\qoq{\quad\text{or}\quad}

\begin{document}

\title[On  schurian fusions of the  scheme of a Galois  plane]{On  schurian fusions of the association scheme of a Galois  affine plane of prime order}
\author{Bahareh Asadian}
\address{Department of pure Mathematics, Faculty of Mathematical Sciences, Shahrekord University, Shahrekord, Iran}
\email{asadian.bahare@gmail.com}
\author{Ilia Ponomarenko}
\address{St.Petersburg Department of the Steklov Mathematical Institute, St.Petersburg, Russia}
\email{inp@pdmi.ras.ru}
\thanks{The work of the second author was supported by the RAS Program of Fundamental Research ``Modern Problems of Theoretical Mathematics''. }

\date{}

\begin{abstract}
The  schurian fusions of the association scheme of a Galois  affine plane of prime order
are completely identified.
\end{abstract}

\maketitle

\section{Introduction}
An {\it association scheme} $\cX$ on a (finite) set $\Omega$ can be thought as a special partition~$S$ of the Cartesian square $\Omega^2$, that contains a diagonal as one of the classes (for the exact definitions, see Section~\ref{290918a}). It is very rare that each  coarser partition of~$\Omega^2$  with the diagonal as a class is also an association scheme, a {\it fusion} of~$\cX$. In~\cite{golfand}, it was proved that this is true if $\cX$ is the  scheme of a finite affine plane~$\cA$, i.e., $\Omega$ is the point set of~$\cA$ and the nondiagonal classes of~$S$ are in one-to-one correspondence with the parallel classes of~$\cA$. Thus if $\cA$ is of order ~$q$, then $|\Omega|=q^2$ and $\cX$ has exactly $p(n)$ different fusions, where $n=q+1$ and $p(n)$ is the number of all partitions of the set $\{1,\ldots,n\}$.

An association scheme $\cX$ on $\Omega$  is said to be {\it schurian} if there exists a group $K\le\sym(\Omega)$ such that the classes of the partition~$S$  are the orbits of the induced action of~$K$ on~$\Omega^2$. The schurity problem in a class of association schemes consists in identifying the schurian schemes in the class in question, see~\cite{ilia}. In the present paper, we solve this problem for the class of all  schurian fusions of the association scheme of a Galois  affine plane of prime order.\medskip

{\bf Main Theorem.} {\it 
A schurian fusion of the  scheme of a Galois affine plane of prime order~$p$ is one of the following:
\nmrt
\tm{1}   wreath or subtensor product of two trivial 	schemes of degree $p$,
\tm{2}	primitive pseudocyclic scheme,
\tm{3}  one of the two exceptional 	schemes,
\tm{4}  the involutive fusion of one of the above schemes.
\enmrt
}\medskip
 
The first three cases in the Main Theorem are basic. In case~(1), the wreath product is unique and schurian, whereas there are non-schurian subtensor products, see example in~\cite[Theorem~26.4]{problem}. The schurian schemes in case~(2) are obtained from $3/2$-transitive subgroups of $\AGL(2,p)$; again there are many non-schurian primitive pseudocyclic schemes, see~\cite[Example~2.6.15]{CP}. Two exceptional schurian schemes from case~(3) correspond to the alternating subgroups~$\alt(4)$ and~$\alt(5)$ of the group $\PGL(2,p)$. For certain values of~$p$, these schemes may be primitive pseudocyclic, see Subsection~\ref{300918a}.

A fusion of a scheme $\cX$ is said to be {\it involutive} if there exists an algebraic automorphism $\varphi$ of $\cX$ such that each class of the partition associated with this fusion is of the form $s\cup\varphi(s)$, $s\in S$. The class of schemes in case~(4) is quite large and can contain schemes occurring in the other three cases. Moreover, many involutive fusions of (even schurian) schemes are non-schurian.
 
The proof of the Main Theorem is given in Sec.~\ref{300918b}; the key ingredients are a classification of $2$-closed permutation groups of prime-squared degree~\cite{dobson}  and an information on the orbits of subgroups of $\PGL(2,q)$ \cite{pgl}.  In Sec.~\ref{290918a}, we cite some standard facts on  association schemes. The  scheme of an affine plane is defined and studied in Sec.~\ref{300918u}. Section~\ref{300918o} contains concluding remarks and open problems.\medskip
 
{\bf Notation.}

Throughout this paper, $\Omega$ is a finite set.

The diagonal of the Cartesian product $\Omega^2$ is denoted by $1_{\Omega}$.  For a relation $s\subseteq \Omega^2$, we set $s^\ast=\{(\beta,\alpha):(\alpha,\beta)\in S\}$
and $\alpha s=\{\beta\in \Omega:(\alpha,\beta)\in s\}$ for all $\alpha\in \Omega$. For $S\subseteq 2^{\Omega^2}$, we  denote by ${\rm S}^{\cup}$ the set of all unions of the elements of~$S$. We define $ S^{\ast}=\{s^{\ast}:s\in S\}$, $S^\#=S\setminus
\{1_{\Omega}\}$ and $\alpha S=\cup_{s\in S}\alpha s$, where $\alpha \in \Omega$. By $C_{p}$ and $\mF_{q}$, we denote the cyclic group of order $p$ and a finite field of order $q$, respectively. By  $\sym(n)$, $\alt(n)$, and $D_{2n}$, we denote the symmetric and alternating group of degree~$n$, and dihedral group of order $2n$, respectively.

\section{Association schemes}\label{290918a}

In this section, we cite all required concepts on association schemes; the notation, terminology  and  results are taken from~\cite{CP}, see also \cite{ilia}.

\subsection{Definitions.}
Let $\Omega$ be a finite set and $S$ a partition of the Cartesian square~$\Omega^2$.
A pair $\cX=(\Omega,S)$ is called an {\it association scheme} or {\it scheme} on $\Omega$ if the following conditions are satisfied: $1_\Omega\in S$,  $S^*=S$, and given $r,s,t\in S$, the number
$$
c_{rs}^t:=|\alpha r\cap \beta s^*|
$$
does not depend on the choice of $(\alpha,\beta)\in t$. The elements of $\Omega$, $S$, $S^\cup$, and the numbers~$c_{rs}^t$ are called the {\it points}, {\it basis relations}, {\it relations}, and {\it intersection numbers} of~$\cX$, respectively. The numbers $|\Omega|$ and $|S|$ are called the {\it degree} and {\it rank} of~$\cX$. A scheme of rank~$2$ is said to be {\it trivial}. The set $S$ of all basis relations of $\cX$ is denoted by $S(\cX)$.

\subsection{Isomorphisms and schurity.}
A bijection from the point set of a scheme $\cX$ to the point set of a scheme~$\cX'$ is
called an {\it isomorphism} from $\cX$ to $\cX'$ if it induces a bijection between their sets of basis relations. The schemes $\cX$ and $\cX'$ are said to be isomorphic
if there exists an isomorphism from $\cX$ to $\cX'$.

An isomorphism from a scheme $\cX$ to itself is called {\it automorphism} if the induced permutation of the basis relations of $\cX$ is the identity. The set 
$$
\aut(\cX)=\{f\in\sym(\Omega):\ s^f=s\ \,\text{for all}\ \, s\in S\}
$$
of all automorphisms of~$\cX$ is a group  with respect to composition. One can easily see that $\aut(\cX)=\sym(\Omega)$ if and only if the scheme $\cX$ is trivial.

Let $K\le\sym(\Omega)$ be a transitive permutation group, and let $S$ denote the set of orbits in the induced action of~$K$ on~$\Omega^2$. Then,
$$
\inv(K):=(\Omega,S)
$$ 
is a scheme; we say that $\inv(K)$ {\it is  associated} with~$K$. A scheme~$\cX$ on $\Omega$ is said to be {\it schurian} if it is associated with the group~$\aut(\cX)$ (or equivalently with a certain  transitive permutation group on~$\Omega$). 

\subsection{Algebraic isomorphisms and fusions.}
Let $\cX$ and $\cX'$ be schemes. A bijection $\varphi:S\to S',\ r\mapsto r'$ is called an {\it algebraic isomorphism} from~$\cX$ to~$\cX'$ if
\qtnl{f041103p1}
c_{rs}^t=c_{r's'}^{t'},\qquad r,s,t\in S.
\eqtn
Each isomorphism~$f$ from~$\cX$ onto~$\cX'$ induces an algebraic isomorphism $s\mapsto s^f$, but not every algebraic isomorphism is induced by an isomorphism. The group of all algebraic automorphisms of $\cX$ is denoted by $\Aaut(\cX)$.

Let $K\le\Aaut(\cX)$. Given $s\in S$, denote by $s^K$ the union of all relations $s^k$, $k\in K$. Then the pair
$$
\cX^K=(\Omega,S^K)
$$
with $S^K=\{s^K:\ s\in S\}$, is called the {\it algebraic fusion} of $\cX$ with respect to the group~$K$. When the order of~$K$ equals~$2$, the fusion is said to be~\emph{involutive}.

\subsection{Parabolics.}
Let $\cX=(\Omega,S)$ be a scheme. Following \cite{Higman1995}, any equivalence relation $e\in S^\cup$ is called a {\it parabolic} of $\cX$. Clearly, $1_\Omega$ and $\Omega^2$ are parabolics of~$\cX$; they are said to be trivial. The scheme~$\cX$ is said to be {\it primitive} if they are the only parabolics of~$\cX$; otherwise, $\cX$ is said to be {\it imprimitive}. The following almost obvious statement is well known.

\prpstnl{200918a}
For a transitive group~$K$, the scheme $\inv(K)$ is primitive if and only if so is the group~$K$.
\eprpstn 

Let $e$ be a parabolic of $\cX$. Denote by $\Omega/e$ the set of all classes of~$e$. For any  $s\in S$, we define $s_{\Omega/e}$  to be the relation on $\Omega/e$ that consists of all pairs $(\Delta,\Gamma)$ such that the relation $s_{\Delta,\Gamma}=s\cap(\Delta\times\Gamma)$
is not empty.  Then the pairs
$$
\cX_{\Omega/e}=(\Omega/e,S_{\Omega/e}) \qaq\cX_\Delta=(\Delta,S_\Delta),
$$
where $S_{\Omega/e}$ and $S_\Delta$ are the sets of all nonempty relations of the form $s_{\Omega/e}$ and~$s_{\Delta,\Delta}$, respectively, are  schemes; here, $s$ runs over $S$, and $\Delta\in\Omega/e$ is fixed.  

If $\cX$ is schurian, then $\cX_{\Omega/e}$  is the scheme associated with the group induced by the action of $\aut(\cX)$ on $\Omega/e$, whereas $\cX_\Delta$ is the scheme induced by the action of the setwise stabilizer of $\Delta$ in $\aut(\cX)$ on~$\Delta$.
 
 \subsection{Wreath and subtensor products.}
 
Let $\Omega_1$ and $\Omega_2$ be sets and  $\Omega=\Omega_1\times\Omega_2$.  Denote by $e_1$ and $e_2$ the equivalence relations on $\Omega$ such that
$$
\Omega/e_1=\{\{\alpha\}\times\Omega_2:\ \alpha\in\Omega_1\}\qaq
\Omega/e_2=\{\Omega_1\times\{\alpha\}:\ \alpha\in\Omega_2\}.
$$
In what follows, the set $\Omega_i$ is canonically identified both with $\Omega/e_i$ and with a class of the equivalence relation $e_{3-i}$, $i=1,2$.
 
Let $\cX_1$ and $\cX_2$ be schemes on $\Omega_1$ and $\Omega_2$, respectively.  The {\it wreath product} of~$\cX_1$ and $\cX_2$ is defined to be the scheme on $\Omega$ that has the smallest rank among the schemes~$\cX$ having a parabolic $e=e_2$ and such that 
$$
\cX_{\Omega_1}=\cX_1\qaq\cX_{\Omega/e_2}=\cX_2,
$$
where $\Omega_1$ on the left-hand side is treated as a class of $e$ (in particular, $\cX_\Delta=\cX_1$ for all $\Delta\in\Omega/e_1$). The basis relations of the wreath product can be found explicitly, see~\cite[Subsection~3.4.1]{CP}.

A {\it subtensor product} of $\cX_1$ and $\cX_2$ is defined to be a scheme $\cX=(\Omega,S)$ such that $e_1$ and $e_2$ are parabolics of~$\cX$,
$$
\cX_{\Omega/e_1}=\cX_1\qaq \cX_{\Omega/e_2}=\cX_2,
$$
and each relation of $\cX$ is contained in the product
$$
s_1\otimes s_2=\{((\alpha_1,\alpha_2),(\beta_1,\beta_2))\in\Omega\times\Omega:\ (\alpha_1,\alpha_2)\in s_1,\ (\beta_1,\beta_2)\in s_2\}, 
$$
where $s_1$ and $s_2$ are basis relations of $\cX_1$ and $\cX_2$, respectively. Such a scheme is not unique and coincides with the tensor product of $\cX_1$ and $\cX_2$ if the rank of $\cX$ equals the product of the ranks of $\cX_1$ and $\cX_2$,  see~\cite[Subsection~3.2.2]{CP}.

\prpstnl{210918a}
Let $K_1\le\sym(\Omega_1)$ and $K_2\le\sym(\Omega_2)$. Then
\nmrt
\tm{1} the scheme of the wreath product $K_1\wr K_2$ in the imprimitive action equals the  wreath product of $\inv(K_1)$ and $\inv(K_2)$, 
\tm{2} the scheme of the subdirect product $K_1\sqcup K_2$ in the product action  equals the  subtensor product of $\inv(K_1)$ and $\inv(K_2)$.
\enmrt
 \eprpstn
\proof Follows from~\cite[Theorem~3.4.6]{CP} and ~\cite[Subsection~3.2.21]{CP}.\eprf

\subsection{Pseudocyclic schemes.} Let $\cX=(\Omega,S)$ be a scheme, and let $s$ be a basis relation of~$\cX$. The numbers
$$
n_s=c_{ss^*}^{1_\Omega}\qaq c(s)=\sum_{r\in S}c_{rr^*}^s
$$
are called the {\it valency} and {\it indistinguishing number} of $s$, respectively. The scheme~$\cX$ is said to be {\it pseudocyclic} if there exists a positive integer~$k$ such that
$$
n_s=k=c(s)+1
$$
for all $s\in S^\#$ (another but equivalent definition is given in \cite[Theorem~3.2]{ilia2}). It is known that the scheme of any Frobenius group is pseudocyclic, and the converse statement is true whenever $|\Omega|$ is much greater than~$k$.

\section{Affine schemes and their fusions}\label{300918u}

Let $\cA$ be a finite affine plane with point set $\Omega$. Then the set $\Omega^2\setminus 1_\Omega$ can be partitioned into the classes according to parallelism: two pairs $(\alpha,\beta)$ and $(\alpha',\beta')$ of points are in one class if
and only if
$$
\alpha\beta=\alpha'\beta'\qoq \alpha\beta\parallel\alpha'\beta',
$$
where $\alpha\beta$ and $\alpha'\beta'$ are the lines through $\alpha$ and $\beta$, and $\alpha'$ and $\beta'$, respectively. 

The obtained classes together with $1_\Omega$ form a partition of $\Omega^2$; denote it by~$S_\cA$. Then the pair 
$$
\cX_\cA=(\Omega,S_\cA)
$$
is an association scheme~\cite{golfand}. It is called the scheme of~$\cA$~\cite{golfand}. The basic properties of this scheme are straightforward and given in the lemma below, see also~\cite{golfand,PR}.

\lmml{060918a}
In the above notation, let $q$ be the order of~$\cA$, $\cX=\cX_\cA$, and $S=S_\cA$. The following statements hold:
\nmrt
\tm{1} $|\Omega|=q^2$ and $|S^\#|=q+1$,
\tm{2}  any $s\in S^\#$ is the disjoint union of $q$ complete graphs of order~$q$; in particular, $n_s=q-1$,
\tm{3} $\Aaut(\cX)=\sym(S)_{1_\Omega}$;\footnote{Here, $\sym(S)_{1_\Omega}$ is the point stabilizer of~$1_\Omega$ in~$\sym(S)$.} in particular, the scheme $\cX$ is pseudocyclic.
\enmrt
\elmm

\crllrl{300918p}
Let $\cX$ be a fusion of the scheme $\cX_\cA$. Then given a parabolic $e$ of~$\cX$ and $\Delta\in\Omega/e$, the schemes $\cX_\Delta$ and $\cX_{\Omega/e}$ are trivial.
\ecrllr

Let $\cX$ be a fusion of the scheme $\cX_\cA$. From statement~(2)  of Lemma~\ref{060918a}, it follows that the valency of any irreflexive basis relation of~$\cX$ is a multiple of $q-1$. Set
$$
\Lambda(\cX)=\Bigl\{\frac{n_s}{q-1}:\ s\in S(\cX)^\#\Bigr\}.
$$
Clearly, this set contains at most $q+1$ positive integers each of which is less than or equal to~$q+1$. 

\lmml{imprim}
In the above notation, set  $\Lambda=\Lambda(\cX)$. Then 
\nmrt
\tm{1} $\cX$ is imprimitive if and only if $1\in \Lambda$,
\tm{2} $\cX$ is pseudocyclic if and only if $|\Lambda|=1$.
\enmrt
\elmm
\proof 	The ``if'' part of statement~(1) immediately follows from statement~(2) of Lemma~\ref{060918a}. To prove the ``only if'' part, assume that the scheme $\cX$ is imprimitive. Then there is a nontrivial parabolic~$e$ of~$\cX$. Denote by~$a$ the number of irreflexive basis relations of $\cX$ contained in $e$. By statements~(1) and~(2) of Lemma~\ref{060918a}, we have
$$
1\le a<q+1\qaq 1+a(q-1)\ \text{divides}\ q^2.
$$
Consequently, $a=1$. It follows that $e=1_\Omega\cup s$ for some $s\in S(\cX)^\#$. Thus, $\Lambda(\cX)$ contains the number $\frac{n_s}{q-1}=1$.

The ``only if'' part of statement~(2) immediately follows from the definition of pseudocyclic scheme. To prove the ``if'' part, assume that $\Lambda(\cX)=\{d\}$ for some positive integer $d\le q+1$. Then each irreflexive basis relation of~$\cX$ is a union of exactly $d$ relations belonging $S_\cA^\#$. By statement~(3) of Lemma~\ref{060918a}, this implies that there exists a cyclic group 
$$
K\le\Aaut(\cX_\cA)
$$ 
of order $d$ that fixes $1_\Omega$, acts semiregularly on $S_\cA^\#$.  Thus in accordance with \cite[Theorem 3.4]{ilia2}, the scheme $\cX$ is pseudocyclic.\eprf\medskip

Let $\cA$ be a Galois affine plane of order~$q$. It is easily seen that the group $\aut(\cX_\cA)$ contains the center of $\GL(2,q)$. Now if $\cX$ is a fusion of $\cX_\cA$, then $\aut(\cX)$ contains $\aut(\cX_\cA)$, and hence
\qtnl{180918d}
Z(\GL(2,q))\le\aut(\cX).
\eqtn

From now on assume that $\cX$ is schurian and, in addition,
\qtnl{180918a}
\aut(\cX)\le\AGL(2,p).
\eqtn  
Then the group $\aut(\cX)$ preserves the parallelism in~$\cA$ and hence acts on the parallel classes of~$\cA$. Since the parallel classes are in one-to-one correspondence with the relations of~$\cS_\cA$, this action induces a group $K\le\sym(\cS_\cA)$ leaving the relation~$1_\Omega$ fixed. By statement~(3) of Lemma~\ref{060918a}, this implies that
$$
K\le\Aaut(\cX_\cA).
$$

Since $K$ is induced by the automorphism group of $\cX$, this scheme is the algebraic fusion of $\cX_\cA$ with respect to~$K$. On the other hand, in view of~\eqref{180918d} and~\eqref{180918a} the group~$K$ can be identified with a subgroup of $\PGL(2,q)$ acting on $q+1$  points of the underlying projective line.
Thus, the following statement holds.

\thrml{pro}
Let $\cA$ be  a Galois affine plane  of order~$q$ and $\cX$ a schurian fusion of~$\cX_\cA$. Assume that condition~\eqref{180918a} holds. Then there is  $K\le \PGL(2,q)$ such that
$$
\cX=(\cX_\cA)^K.
$$ 
In particular, $\Lambda(\cX)$ equals the set $N(K)$ of cardinalities of the orbits of~ $K$.
\ethrm

\section{ The proof of the Main Theorem}\label{300918b}

By the hypothesis of the theorem, $\cX$ is the scheme of the group $\aut(\cX)$; in particular, $\cX$ is primitive (respectively, imprimitive) if and only if $\aut(\cX)$  is primitive (respectively, imprimitive) (Proposition~\ref{200918a}). The proof is divided into two parts depending on whether or not the group scheme~$\cX$ is imprimitive. 

The imprimitive case corresponds to statement~(1) of the Main Theorem; here we use  a characterization of the $2$-closed subgroups of $\sym(p^2)$  given in~\cite{dobson}. Statements~(2), (3), and~(4) of the Main Theorem arise in the primitive case; here our tool is the information on the subgroups of $\PGL(2,q)$ given in~\cite{pgl}.

\subsection{The scheme $\cX$ is imprimitive.} The group~$\aut(\cX)$ being the automorphism group of a scheme is $2$-closed in the sense of \cite{Wielandt1969}. Therefore, we make use of the following statement which is an immediate consequence of~\cite[Theorem 14]{dobson}.

\lmml{doblem}
Let $K\le \sym(p^2)$ be a $2$-closed  group with a regular subgroup~$C_p\times C_p$. Then one of the following statements holds.
\nmrt
\tm{i} $K$ is primitive, and $K\leq \AGL(2, p)$,  or $K=\sym(p)\wr \sym(2)$ or  $\sym(p^2)$,
\tm{ii} $K$ is imprimitive, and one of the following statements holds:
\nmrt
\tm{ii1} $K=\sym(p)\times K'$, where $K'\le\sym(p)$,
\tm{ii2} $K<\AGL(1, p)\times \AGL(1, p)$,
\tm{ii3} $K=K_1\wr K_2$, where $K_1,K_2\le\sym(p)$ are $2$-closed groups.
\enmrt
\enmrt
\elmm

By Lemma~\ref{doblem} for $K=\aut(\cX)$, we have two cases: the first one is formed by  statements~(ii1) and~(ii2), whereas the second one consists of just  statement~(ii3). In the former case, $K$  is subdirect product of two groups. Therefore the scheme~$\cX$ is the subtensor product of  two  schemes of degree $p$ (statement~(2) of Proposition~\ref{210918a}), and both of them are trivial (Corollary~\ref{300918p}). In the latter case, $\cX$ is  the wreath product $\inv(K_1)\wr\inv(K_2)$ (statement~(1) of Proposition~\ref{210918a}), and again both of them are trivial (Corollary~\ref{300918p}). Thus if $\cX$ is imprimitive, then statement~(1) of the Main Theorem holds.

\subsection{The scheme $\cX$ is primitive.} Without loss of generality, we may assume that (a) $\cX$ is not trivial, for otherwise statement~(2) of the Main Theorem holds and (b) the relation 
\qtnl{010918a}
1\not\in\Lambda(\cX)
\eqtn
holds, for otherwise $\cX$ is imprimitive by statement~(1) of Lemma~\ref{imprim}. Then $p$ is odd and the following statement is a special case of the results proved in ~\cite[Theorem~2 and  Sec.~4]{pgl}. 

\lmml{020918a}
Let $K\le\PGL(2,p)$ be an intransitive permutation group acting on $p+1$  points of the underlying projective line, and $N=N(K)$. Then one of the following statements holds:
\nmrt
\tm{1} $K=C_d$ and $N\subseteq\{1,d\}$, $d\ge 1$, 
\tm{2} $K=D_{2d}$ and $N\subseteq\{2,d,2d\}$, $d\geq 2$,  
\tm{3}$K=C_p\rtimes C_d$ and $N\subseteq\{1,p\}$, $d\mid p-1$, 
\tm{4} $K=\alt(4)$,  $\alt(5)$, or $\sym(4)$.
\enmrt
\elmm

By Theorem~\ref{pro} for $q=p$, there exists a group $K$ satisfying the hypothesis of Lemma~\ref{020918a} and such that  
$$
\cX=(\cX_\cA)^K\qaq\Lambda=N,
$$
where $\cA$ is a  Galois affine plane of order~$p$ and $\Lambda=\Lambda(\cX)$. Note that this group is intransitive, because the scheme $\cX$ is nontrivial. To complete the proof we will verify that in each of the four cases of Lemma~\ref{020918a}, the conclusion of the Main Theorem holds.

In the case~(1), assumption~\eqref{010918a} implies that $N=\{d\}$. It follows that $|\Lambda|=1$. Thus the scheme $\cX$ is pseudocyclic by statement~(2) of Lemma~\ref{imprim}.

In the case~(2), one can see as above that the scheme $\cX$ is pseudocyclic whenever $2\not\in N$ and $d\not\in N$. Assume first that $2\in N$. Denote by $K'$ the kernel of the action of~$K$ on an orbit of size~$2$.  Then $K'$ is a subgroup of index~$2$ and  $1\in N(K')$. It follows  that if 
\qtnl{020918u}
\cX'=(\cX_\cA)^{K'},
\eqtn
then $\cX$ is an involutive  fusion of $\cX'$ and  $1\in N(K')= \Lambda(\cX')$. The scheme $\cX'$ is imprimitive by statement~(1) of Lemma~\ref{imprim}. By the first part of the proof (the imprimitive case), this implies that statement~(1) of the Main Theorem holds for~$\cX'$, and we are done.

Remaining in the case~(2), we may assume that $N=\{d,2d\}$. Then $K$ has a subgroup $K'$ of index~$2$ such that
\qtnl{020918f}
N(K')=\{d\}.
\eqtn
Indeed, the action of $K$ on an orbit of cardinality $d$ is permutation isomorphic to the action of $K$ on the right cosets of a subgroup generated by an involution $k\in K$.  Depending on whether or not $k$ lies in the center of~$K$, one can take as~$K'$ a subgroup of $K$ isomorphic to $D_d$ or~$C_d$.  Now, in view of~\eqref{020918f}, the scheme $\cX'$ defined by formula~\eqref{020918u} is pseudocyclic (statement~(2) of Lemma~\ref{imprim}). Therefore statement~(2) of the Main Theorem holds for~$\cX'$. Since $\cX$ is an involutive fusion of~$\cX'$, we are done. 

To complete the proof, it suffices to note that in the case~(3) the scheme $\cX$ is pseudocyclic by assumption~\eqref{010918a}, whereas in the case~(4) the scheme $\cX$ is either exceptional ($K=\alt(4)$ or $\alt(5)$), or an involutive fusion of the scheme~\eqref{020918u} with $K'=\alt(4)$ for $K=\sym(4)$.

\section{Concluding remarks}\label{300918o}

In what follows,  $\cC_1$, $\cC_2$, $\cC_3$, and $\cC_4$   denote the classes  of  schemes in statements~(1), (2), (3), and~(4) of the Main Theorem, respectively. 

\subsection{Interrelation between the classes from the Main Theorem.}\label{300918a} In view of  the remarks made after the Main Theorem, we are interested in the interrelation between the classes $\cC_1$, $\cC_2$, and $\cC_3$. The schemes in $\cC_1$ are imprimitive, whereas those in $\cC_2$ and~$\cC_3$ are not. Therefore,
$$
\cC_1\cap\cC_2=\cC_1\cap\cC_3=\varnothing.
$$
The classes $\cC_2$ and $\cC_3$ have nontrivial intersection. This follows from the information on the orbit lengths of the groups $\alt(4),\alt(5)\le \PGL(2,p)$ obtained in~\cite[Lemmas~9,11]{pgl}. Indeed, the exceptional schemes associated with groups $\alt(4)$ and $\alt(5)$ are primitive pseudocyclic if, e.g., 
$$
\mmod{p}{-1}{a}, \quad a=3\ \,\text{or}\ \,5.
$$

\subsection{The automorphism groups.}	In principle, all the information of the automorphism group of the scheme $\cX$ in the Main Theorem can be extracted from Lemma~\ref{doblem}. In the most cases, we have
$$
\aut(\cX)\le\AGL(2,p),
$$
i.e., $\cX$ is isomorphic to a normal Cayley scheme over $C_p\times C_p$ in the sense of~\cite{circulant1}. Apart from this case, the only possibility for the group $\aut(\cX)$ are the following:
\qtnl{030918i}
\sym(p)\times\sym(p),\quad\sym(p)\wr\sym(p),\quad \sym(p)\wr\sym(2),\quad\sym(p^2).
\eqtn
The first two groups appear in statements (ii1) and (ii3) of Lemma~\ref{doblem} and the schemes of these groups are in the class~$\cC_1$, whereas the second two groups appear in statement~(i) and the schemes of these groups are the Hamming scheme $H(2,p)$ and trivial scheme lying in the classes $\cC_4$ and $\cC_2$, respectively.

\subsection{Further research.} The first natural problem is to generalize the Main Theorem to the $p$-powers~$q$, i.e., to find a compact description of schurian fusions of a Galois affine plane of order~$q$. In this way, one can still use the results of~\cite{pgl} where  they were established  arbitrary~$q$. However, to the author knowledge, there is no  generalization of Lemma~\ref{doblem}.

The class $\cC_2$ contains the cyclotomic schemes over near-fields of order~$p^2$ \cite{BP08}  and the schemes of Frobenius groups. It would be interesting to find other schemes in~$\cC_2$ (if they are).

From the algorithmic point of view, one of the problem in the above context is how to recognize the schemes $\cX$ from the Main Theorem  in the class of all association schemes efficiently. Definitely, this can easily be done if $\aut(\cX)$ is one of the groups~\eqref{030918i}. For the other schemes, the problem can efficiently be reduced to recognizing schemes belonging to the classes~$\cC_2$ and~$\cC_4$.


\begin{thebibliography}{99}


\bibitem{BP08} 
J.~Bagherian, I.~Ponomarenko,   and A.~Rahnamai~Barghi, ``On cyclotomic schemes over finite near-fields,'' {\it J. Algebraic Combin.}, {\bf 27}, No.~2, 173--185 (2008).
	

\bibitem{pgl} 
P.~J.~Cameron, G.~R.~Omidi, and B.~Tayfeh-Rezaie, ``$3$-designs from $\PGL(2,q)$,'' \textit{ Electronic J. Combin.}, \textbf{13}, No.~1, \#50 (2006).
 
\bibitem{CP}
G.~Chen and I.~Ponomarenko, \emph{Lectures on Coherent Configurations} (2019),
{\tt http://www.\-pdmi.ras.ru/\textasciitilde inp/ccNOTES.pdf}.
 
\bibitem{dobson} 
E. Dobson and D. Witte, ``Transitive permutation groups of prime-squared degree,'' \textit{J. Algebraic Combin.}, {\bf 16}, No.~1, 43--69 (2002).

\bibitem{circulant1} 
S.~Evdokimov and I.~Ponomarenko, ``Characterization of cyclotomic schemes and normal Schur rings over a cyclic group,''  {\it St. Petersburg Math J.}, \textbf{14}, No.~2, 189--221 (2002).


\bibitem{ilia} 
S.~Evdokimov and I.~Ponomarenko, ``Permutation group approach to association schemes,'' \textit{Eroupean J. Combin.}, {\bf 30}, No.~6, 1456--1476 (2009).

\bibitem{golfand}
A.~Yu.~Gol'fand, A.~V.~Ivanov, and M.~Kh.~Klin, ``Amorphic cellular rings,'' in: {\it Investigations in algebraic theory of combinatorial objects}, Kluwer
Acad. Publ., Dordrecht (1994), pp.~167--186.

\bibitem{Higman1995}
D.~G.~Higman, ``Rank $5$ association schemes and triality,'' {\it Linear Algebra Appl.}, \textbf{226/228}, 197--222 (1995).
 

\bibitem{ilia2} 
M. Muzychuk and I. Ponomarenko, ``On pseudocyclic association schemes,'' \textit{ARS Math.  Contemp.}, {\bf 5}, 1--25 (2012).

\bibitem{PR} 
I.~Ponomarenko and  A.~Rahnamai~Barghi, ``On amorphic C-algebras,'' {\it J. Math. Sci.}, {\bf 145},  No.~3, 4981--4988  (2007).

\bibitem{problem} 
H.~Wielandt, {\it Finite Permutation Groups}, Academic Press, New York- London (1964).

\bibitem{Wielandt1969}
H.~Wielandt, {\it Permutation Groups Through Invariant Relations and Invariant 	Functions}, The Ohio State University (1969).

\end{thebibliography}
\end{document}